\documentclass[12pt]{amsart}
\usepackage{amssymb}
\usepackage[all]{xy}

\DeclareMathOperator{\Hom}{Hom}
\DeclareMathOperator{\End}{End}
\DeclareMathOperator{\Aut}{Aut}
\DeclareMathOperator{\SL}{SL}
\DeclareMathOperator{\GL}{GL}
\DeclareMathOperator{\PGL}{PGL}
\DeclareMathOperator{\Sp}{Sp}
\DeclareMathOperator{\Out}{Out}
\DeclareMathOperator{\Pic}{Pic}

\DeclareMathOperator{\ad}{ad}
\DeclareMathOperator{\codim}{codim}
\DeclareMathOperator{\Spec}{Spec}
\DeclareMathOperator{\Sym}{Sym}

\newtheorem{theorem}{Theorem}[section]
\newtheorem{proposition}[theorem]{Proposition}
\newtheorem{lemma}[theorem]{Lemma}
\newtheorem{corollary}[theorem]{Corollary}
\theoremstyle{remark}
\newtheorem{remark}[theorem]{Remark}

\newcommand{\longto}[1][]{\stackrel{#1}{\longrightarrow}}
\newcommand{\bbZ}{\mathbb{Z}}
\newcommand{\bbC}{\mathbb{C}}
\newcommand{\bbP}{\mathbb{P}}
\newcommand{\Gm}{\bbC^*}
\newcommand{\calM}{\mathcal M}
\newcommand{\calL}{\mathcal L}
\newcommand{\calO}{\mathcal O}
\newcommand{\pr}{\mathrm{pr}}
\newcommand{\sing}{\mathrm{sing}}
\newcommand{\stab}{\mathrm{stab}}

\newcommand{\Gbar}{\overline{G}}
\newcommand{\Hbar}{\overline{H}}
\newcommand{\Ubar}{\overline{U}}
\newcommand{\ibar}{\overline{\textit{\i}}}
\newcommand{\dbar}{\bar{d}}
\newcommand{\ebar}{\bar{e}}

\newcommand{\prbar}{\overline{\pr}}
\newcommand{\Xtilde}{\widetilde{X}}
\newcommand{\Stilde}{\widetilde{S}}
\newcommand{\Ldet}{\mathcal{L}_{\det}}

\numberwithin{equation}{section}

\begin{document}

\title[Torelli for moduli of principal bundles]{A Torelli theorem for
moduli spaces of principal bundles over a curve}

\author[I. Biswas]{Indranil Biswas}

\address{School of Mathematics, Tata Institute of Fundamental
Research, Homi Bhabha Road, Bombay 400005, India}

\email{indranil@math.tifr.res.in}

\author[N. Hoffmann]{Norbert Hoffmann}

\address{Freie Universit\"at Berlin, Institut f\"ur Mathematik,
Arnimallee 3, 14195 Berlin, Germany}

\email{norbert.hoffmann@fu-berlin.de}



\date{}

\begin{abstract}
Let $X$ and $X'$ be compact Riemann surfaces of genus $\geq 3$,
and let $G$ and $G'$ be nonabelian reductive complex groups.
If one component $\calM_G^d( X)$ of the coarse moduli space for
semistable principal $G$--bundles over $X$ is isomorphic to
another component $\calM_{G'}^{d'}(X')$, then $X$ is isomorphic
to $X'$.
\end{abstract}

\maketitle

\section*{Introduction}
Let $X$ be a compact connected Riemann surface of genus $g_X$.
The classical Torelli theorem says that the isomorphism
class of $X$ is uni\-que\-ly determined by the isomorphism class
of the Jacobian $\Pic^0( X)$, together with
its canonical principal polarization $\Theta$.

Several authors have studied nonabelian analogues, replacing
Jacobians by moduli spaces of vector bundles. Suppose $g_X \geq 3$.
Given a line bundle $L$ on $X$, let $\calM_{n, L}( X)$
denote the moduli space of semistable vector bundles
$E$ over $X$ of rank $n$ with fixed determinant $\bigwedge^n E \cong L$.
Then the isomorphism class of $X$ is uniquely determined by
the isomorphism class of the projective variety $\calM_{n, L}( X)$.

This was first proved for $n = 2$ and $\deg( L)$ odd, 
by Mumford-Newstead \cite[p. 1201, Corollary]{MN} and Tyurin \cite[Theorem 1]{T1}.

It was then proved more generally for $n$ coprime to $\deg( L)$,
by Tyurin \cite[Theorem 1]{T2} and Narasimhan-Ramanan \cite[Theorem 3]{NR_Def}.

Finally, this statement was proved in full generality by Kouvidakis-Pantev \cite[{Theorem E}]{KP},
and then later for $g_X \geq 4$ also by Hwang-Ramanan \cite[Theorem 5.1]{HR} and by Sun \cite[Corollary 1.3]{Su}.

These proofs used either an intermediate Jacobian of the moduli space 
\cite{MN, NR_Def}, or
Higgs bundles and the Hitchin map \cite{KP, HR}, or rational curves on the moduli spaces \cite{T1, T2, Su}.

Our aim here is to address a similar question for the moduli spaces of principal bundles.
Given a connected reductive complex linear algebraic group $G$, let $\calM_G^d( X)$ denote
the moduli space of semistable principal $G$--bundles over $X$ of topological type $d \in \pi_1( G)$.
\begin{theorem} \label{mainthm}
Let $G$ and $G'$ be nonabelian connected reductive complex groups.
Let $d \in \pi_1( G)$ and $d' \in \pi_1( G')$ be given.
Let $X$ and $X'$ be compact Riemann surfaces of genus $\geq 3$.
If $\calM_G^d( X)$ is isomorphic to $\calM_{G'}^{d'}( X')$, then $X$ is isomorphic to $X'$.
\end{theorem}

For the proof, our strategy is fundamentally different from the earlier ones.
Whereas \cite{KP} and \cite{HR} start from the stable locus,
our starting point is the strictly semistable locus.
We show that it lies in the singular locus and
is characterized by its type of singularities.
Using powers of its anticanonical line bundle, we map it to a projective space. 
The fibers of this map allow us to reconstruct the Jacobian of $X$ and its principal polarization.
Then the usual Torelli theorem applies.

This method does not work if the strictly semistable locus is empty. But this case is rare,
and can in fact be reduced to the earlier results on moduli spaces of vector bundles in the coprime case.

The structure of this paper is as follows. Section \ref{sec:singularities} deals with quotient singularities
arising from a group action. We give a criterion to distinguish the case of a finite group from the torus $\Gm$.

In Section \ref{sec:moduli}, we collect some basic facts about the moduli spaces $\calM_G^d( X)$.
We use several tools valid only in characteristic $0$, like representations of $\pi_1( X)$,
and the stability of induced bundles.

Section \ref{sec:str_semistab} describes the strictly semistable locus in $\calM_G^d( X)$.
We relate it to moduli spaces of principal $H$--bundles for Levi subgroups $H$ of maximal parabolic subgroups in $G$.
Furthermore, we characterize it by the singularities it contains. As a byproduct, we show that the smooth locus
of $\calM_G^d( X)$ coincides with the regularly stable locus.

In the final Section \ref{sec:proof}, we prove Theorem \ref{mainthm}, by reconstructing the compact Riemann surface $X$
from the projective variety $\calM_G^d( X)$. In fact we prove a slightly stronger result, namely that $X$
is determined up to isomorphy by the smooth locus of $\calM_G^d( X)$; see Theorem \ref{thm:main}.

\section*{Acknowledgements}
We are very grateful to the referee for pointing out a gap in an earlier version of this paper,
and also for some valuable suggestions. We thank D. Prasad for a useful conversation.
The first author wishes to thank the Freie Universit\"at Berlin and the Institute for Mathematical Sciences at Chennai for hospitality.
The second author was supported by the SFB 647: Raum - Zeit - Materie.
He thanks TIFR (Mumbai) and IMSc (Chennai) for hospitality while some part of this work was done.

\section{Some quotient singularities} \label{sec:singularities}
Let $S$ and $S'$ be schemes of finite type over $\bbC$. We say that $S$ is near a closed point $s_0 \in S$ \emph{analytically isomorphic}
to $S'$ near a closed point $s_0' \in S_0$ if for some open neighborhoods $s_0 \in U \subseteq S$ and $s_0' \in U' \subseteq S'$ in the Euclidean topology,
there exists a biholomorphic isomorphism $U \to U'$ that maps $s_0$ to $s_0'$.

We always denote by $S^{\sing} \subseteq S$ the singular locus of $S$. Let $G$ be a reductive linear algebraic group over $\bbC$.
We assume that $G$ acts linearly on a finite-dimensional complex vector space $V$. Let
\begin{equation*}
  p: V \longto S := V/\!/G = \Spec (\Sym ( V^*)^G)
\end{equation*}
be the GIT--quotient, and put $s_0 := p( 0) \in S$.
\begin{lemma} \label{lemma:finitequotient}
  Let $G$ be finite and nontrivial.
  Suppose that the fixed locus $V^g \subseteq V$ has codimension $\geq 2$ for all $g \in G$ with $g \neq 1$.
  \begin{itemize}
   \item[i)] The quotient $S = V/\!/G$ is singular in $s_0 = p( 0)$.
   \item[ii)] If $s_0 \in U \subseteq S$ is an open neighborhood in the Euclidean topology such that $U \setminus S^{\sing}$ is connected,
    then $\pi_1( U \setminus S^{\sing})$ is nontrivial.
  \end{itemize}
\end{lemma} 
\begin{proof}
  Using the assumption on $\codim V^g$, \cite[Lemma 4.4]{NR} implies
  \begin{equation*}
    p^{-1}( S^{\sing}) = \bigcup_{g \neq 1} V^g.
  \end{equation*}
  This proves (i). In the situation of (ii), $p^{-1}( U \setminus S^{\sing})$ is connected because $s_0$ has only one inverse image.
  Hence $p^{-1}( U \setminus S^{\sing})$ is a nontrivial covering of $U \setminus S^{\sing}$, which is therefore not simply connected.
\end{proof}
\begin{lemma} \label{lemma:Gm}
  Let $G = \Gm$ act linearly on the finite-dimensional complex vector space $V$, with associated weight space decomposition
  \begin{equation*}
    V = \bigoplus_{m \in \bbZ} V_m.
  \end{equation*}
  Suppose that $V_{-1}$ and $V_1$ both have dimension $\geq 2$.
  \begin{itemize}
   \item[i)] The quotient $S = V/\!/\Gm$ is singular in $s_0 = p( 0)$.
   \item[ii)] Every neighborhood $s_0 \in U \subseteq S$ in the Euclidean topology contains an open neighborhood $s_0 \in U' \subseteq U$
    such that $U' \setminus S^{\sing}$ is connected and $\pi_1( U' \setminus S^{\sing})$ is trivial.
  \end{itemize}
\end{lemma}
\begin{proof}
  It is easy to verify $S = (V_+ \oplus V_-)/\!/\Gm \times V_0$ with
  \begin{equation*}
    V_+ := \bigoplus_{m > 0} V_m \qquad\text{and}\qquad V_- := \bigoplus_{m < 0} V_m.
  \end{equation*}
  Replacing $V$ by $V_+ \oplus V_-$ if necessary, we may thus assume $V_0 = 0$ without loss of generality.
  Then the GIT--stable locus $V^{\stab} \subseteq V$ is
  \begin{equation*}
    V^{\stab} = V \setminus ( V_+ \cup V_-),
  \end{equation*}
  and the image $p( V_+ \cup V_-)$ of $V \setminus V^{\stab}$ is just the point $s_0 \in S$. Let
  \begin{equation*}
    V^{\mu_n} = \bigoplus_{m \in n \bbZ} V_m \subseteq V
  \end{equation*}
  be the fixed locus of the finite subgroup $\mu_n \subseteq \Gm$. We claim
  \begin{equation} \label{eq:claim}
    p^{-1}( S^{\sing}) = V_+ \cup V_- \cup \bigcup_{n \geq 2} V^{\mu_n}.
  \end{equation}

  To check this, we first consider a stable point $v \in V^{\stab}$.
  Its image $p( v) \in S$ is smooth if and only if its stabilizer in $\Gm$ is trivial,
  according to Luna's \'{e}tale slice theorem \cite{Lu} and Lemma \ref{lemma:finitequotient}.

  Suppose next that $V^{\mu_n} \cap V^{\stab}$ is non-empty for some $n \geq 2$.
  Then it is dense in $V^{\mu_n}$, but also contained in the closed subset $p^{-1}( S^{\sing})$.
  Hence $0 \in V^{\mu_n} \subseteq p^{-1}( S^{\sing})$, which implies $s_0 \in S^{\sing}$.
  It follows that $p^{-1}( S^{\sing})$ contains the unstable locus $V_+ \cup V_-$ in this case.

  Now suppose $V^{\mu_n} \subseteq V_+ \cup V_-$ for all $n \geq 2$. Then the canonical map
  \begin{equation*}
    \Stilde := (V \setminus \bigcup_{n \geq 2} V^{\mu_n}) \big/ \Gm \longto S
  \end{equation*}
  is surjective. Here $\Gm$ acts freely, so the set $\Stilde$ is a (not necessarily separated) complex manifold.
  Applying \cite[Lemma 4.4]{NR} to this map, we get $s_0 \in S^{\sing}$. Hence $p^{-1}( S^{\sing})$ contains $V_+ \cup V_-$ also in this case.
  This proves the claim \eqref{eq:claim}, and in particular part (i) of the lemma.

  Let a Euclidean neighborhood $s_0 \in U \subseteq S$ be given.
  For each $m \in \bbZ$, we choose a basis $(x_{m, i})_{1 \leq i \leq \dim V_m}$ of the vector space $( V_m)^*$.
  The $\bbC$-algebra $\Sym ( V^*)^{\Gm}$ is generated by finitely many nonconstant monomials $f_1, \ldots, f_N$ in the variables $x_{m, i}$.
  They provide a closed embedding $S \hookrightarrow \bbC^N$ with $s_0 \mapsto 0$. Hence $U$ contains the open neighborhood
  \begin{equation*}
    s_0 \in U' := \{ s \in S: |f_n( s)| < \varepsilon \text{ for } n = 1, \ldots, N\}
  \end{equation*}
  for some $\varepsilon > 0$. Because the $f_n$ are monomials, the open subset
  \begin{equation*}
    p^{-1}( U') = \{v \in V: |f_n( v)| < \varepsilon \text{ for } n = 1, \ldots, N\} \subseteq V
  \end{equation*}
  is a star-shaped neighborhood of $0$, in the sense that $\lambda v \in p^{-1}( U')$ holds for all $v \in p^{-1}( U')$ and
  $\lambda \in [0, 1]$. In particular, $p^{-1}( U')$ is contractible. Since we are only removing finitely many
  linear subspaces of complex codimension $\geq 2$, it follows that $p^{-1}( U' \setminus S^{\sing})$ is connected and simply connected.
  Consequently, $U' \setminus S^{\sing}$ is also connected and simply connected, because
  \begin{equation*}
    p: p^{-1}( U' \setminus S^{\sing}) \longto U' \setminus S^{\sing}
  \end{equation*}
  is a fibration with connected fiber $\Gm$.
\end{proof}

\section{Moduli spaces of principal bundles} \label{sec:moduli}
Let $X$ be a compact connected Riemann surface of genus $g_X \geq 3$.
Let $G$ from now on be a reductive connected linear algebraic group over $\bbC$.
We denote by
\begin{equation*}
  \calM_G = \calM_G( X)
\end{equation*}
the coarse moduli space of semistable algebraic principal $G$--bundles $E$ over $X$. Its connected components
\begin{equation*}
  \calM_G^d = \calM_G^d( X) \subseteq \calM_G( X)
\end{equation*}
are irreducible normal projective varieties of dimension
\begin{equation*}
  \dim \calM_G^d = (g_X-1) \dim_{\bbC} G + \dim_{\bbC} Z_G
\end{equation*}
where $Z_G$ denotes the center of $G$. These connected components are indexed by elements $d \in \pi_1( G)$,
which correspond to topological types of algebraic principal $G$--bundles $E$ over $X$.
In the case $G = \GL_n$ of vector bundles, the topological type $d \in \pi_1( \GL_n) = \bbZ$ is their degree. 

The construction of such moduli spaces is carried out in \cite{Ra_Proc};
see also \cite{GS} and \cite{Sch} for generalizations to higher dimensions.
We recall from \cite{Ra_Proc} that one has a dense open subscheme
\begin{equation*}
  \calM_G^{\stab} \subseteq \calM_G
\end{equation*}
whose closed points correspond bijectively to isomorphism classes of stable principal $G$--bundles over $X$. Closed points
in $\calM_G$ correspond bijectively to isomorphism classes of polystable principal $G$--bundles over $X$, and also to
$S$--equivalence classes of semistable principal $G$--bundles over $X$, in the sense of \cite[Definition 3.6]{Ra_Proc}.

Let $\mathfrak g$ be the Lie algebra of $G$. Given a principal $G$--bundle $E$ over $X$, we denote the associated adjoint vector bundle over $X$ by
\begin{equation*}
  \ad( E) := E \times^G {\mathfrak g} := ( E \times {\mathfrak g})/G.
\end{equation*}
If $E$ is semistable, then $\ad( E)$ is also semistable \cite[Corollary 3.18]{Ra_Proc}. The vector space
\begin{equation*}
  H^1( X, \ad( E))
\end{equation*}
parameterizes infinitesimal deformations of $E$ according to standard deformation theory.
On it, we have an adjoint action of the automorphism group $\Aut( E)$.
If $E$ is polystable, then the GIT--quotient
\begin{equation*}
  H^1( X, \ad( E))/\!/\Aut( E)
\end{equation*}
is near $0$ analytically isomorphic to $\calM_G$ near the point $[E]$. This is a standard consequence of Luna's \'{e}tale slice theorem \cite{Lu}.

The center $Z_G$ of $G$ is a normal subgroup of $\Aut( E)$ for every principal $G$--bundle $E$ over $X$.
$E$ is stable or polystable or semistable if and only if the induced principal  $G/Z_G$--bundle $E/Z_G$ over $X$ is so,
since both bundles have the same set of reductions to parabolic subgroups; cf. \cite[Proposition 7.1]{Ra_Ann}.
If $E$ is stable, then $\Aut( E)/Z_G$ is finite by \cite[Proposition 3.2]{Ra_Ann}.
If moreover $\Aut( E) = Z_G$, then $E$ is called \emph{regularly stable}; in this case, $\calM_G$ is smooth at $[E]$.

The canonical exact sequence of reductive groups
\begin{equation*}
  1 \longto {[G, G]} \longto G \longto G/[ G, G] \longto 1
\end{equation*}
shows that $\pi_1( [G, G])$ is the torsion subgroup of $\pi_1( G)$. We put
\begin{equation*}
  \calM_G' := \bigcup_{d \in \pi_1( [ G, G])} \calM_G^d \subseteq \calM_G.
\end{equation*}
If $G$ is semisimple, then $\calM_G' = \calM_G$.

Choose a maximal compact subgroup $K_G \subseteq G$, and a base-point $x_0 \in X$.
Let $\Xtilde \twoheadrightarrow X$ be the universal covering of $X$.
Given a group homomorphism $\rho: \pi_1( X, x_0) \to K_G$, the principal $G$--bundle
\begin{equation*}
  E_{\rho} := \Xtilde \times^{\rho} G := (\Xtilde \times G)/\pi_1( X, x_0)
\end{equation*}
is polystable, and its moduli point $[E_{\rho}] \in \calM_G$ is contained in $\calM_G'$. In this way, closed points
in $\calM_G'$ correspond bijectively to homomorphism $\rho: \pi_1(X,x_0) \to K_G$ up to conjugacy \cite[Corollary 3.15.1]{Ra_Proc}.
\begin{lemma} \label{lemma:dense}
  The set of points $[E_{\rho}] \in \calM_G'$ coming from homomorphisms $\rho: \pi_1( X, x_0) \to K_G$
  with dense image is Zariski-dense in $\calM_G'$.
\end{lemma}
\begin{proof}
  There are only countably many conjugacy classes of closed subgroups
  $K \subsetneq K_G$, according to \cite[Proposition 10.12]{Ad}.
  For any given closed subgroup $K \subsetneq K_G$, the conjugacy classes of
  homomorphisms $\rho: \pi_1( X, x_0) \to K_G$ factoring through $K$
  form a closed real analytic subset of $\calM_G'$ of a smaller dimension.
  Due to the Baire category theorem, the union of these subsets in
  $\calM_G'$ has empty interior (with respect to the Euclidean topology).
  This shows that the homomorphisms $\rho: \pi_1( X, x_0) \to K_G$
  with dense image are dense in $\calM_G'$ for the Euclidean topology,
  and a fortiori for the Zariski topology.
\end{proof}
\begin{lemma} \label{lemma:Aut}
  Consider the polystable principal $G$--bundle $E_{\rho}$ over $X$ given by a homomorphism $\rho: \pi_1( X, x_0) \to K_G$.
  Its automorphism group $\Aut( E_{\rho})$ is canonically isomorphic to
  \begin{equation*}
    C_G( \rho) := \{g \in G: g \cdot \rho( \omega) = \rho( \omega) \cdot g \text{ for all } \omega \in \pi_1( X, x_0)\}.
  \end{equation*}
\end{lemma}
\begin{proof}
  Every element $g \in C_G( \rho)$ defines a $\pi_1( X, x_0)$--equivariant automorphism of the trivial bundle
  $(\Xtilde \times G)$ over $\Xtilde$, which descends to an automorphism of $E_{\rho}$.
  This embeds $C_G( \rho)$ into $\Aut( E_{\rho})$.

  Now let $\varphi \in \Aut( E_{\rho})$ be given.
  Choose an embedding $G \hookrightarrow \GL_r$ that maps $K_G$ to the unitary group $U_r$.
  Then $\varphi$ defines a nonzero holomorphic section of the unitary vector bundle
  \begin{equation*}
    \Xtilde \times^{\rho} \End( \bbC^r)
  \end{equation*}
  over $X$. This vector bundle is polystable of degree $0$, so the section lies in a trivial direct summand $\calO_X$.
  Such a trivial direct summand corresponds, under the Narasimhan-Seshadri correspondence, to a line in $\End( \bbC^r)$
  on which $\pi_1( X, x_0)$ acts trivially.

  Hence the section in question comes from a fixed point in $\End( \bbC^r)$.
  But this fixed point has to be in the image of $G \hookrightarrow \GL_r \subseteq \End( \bbC^r)$.
  This shows that $\varphi$ indeed comes from an element of $C_G( \rho)$.
\end{proof}
\begin{proposition} \label{prop:fin_quot}
  Let $E$ be a principal $G$--bundle over $X$. Suppose that $E$ is stable, but not regularly stable.
  \begin{itemize}
   \item[i)] The point $[E] \in \calM_G$ is singular.
   \item[ii)] There is a neighborhood $[E] \in U \subseteq \calM_G$ in the Euclidean topology such that
    for every open neighborhood $[E] \in U' \subseteq U$ with $U' \setminus \calM_G^{\sing}$ connected, $\pi_1( U' \setminus \calM_G^{\sing})$ is nontrivial.
  \end{itemize}
\end{proposition}
\begin{proof}
  $\calM_G$ is near $[E]$ analytically isomorphic to the quotient of
  \begin{equation*}
    V := H^1( X, \ad( E))
  \end{equation*}
  modulo the finite group $\Aut( E)/Z_G$ near the image of $0 \in V$.
  Let an element $\varphi \in \Aut( E) \setminus Z_G$ be given.
  Using Lemma \ref{lemma:finitequotient}, it suffices to verify that the fixed locus $V^{\varphi} \subseteq V$ has codimension $\geq 2$.

  The image $K_{G/Z_G} \subseteq G/Z_G$ of $K_G \subseteq G$ is again a maximal compact subgroup.
  Since $E$ is stable, the induced principal $G/Z_G$--bundle $E/Z_G$ is also stable.
  Hence we may assume $E/Z_G = E_{\rho}$ for some group homomorphism $\rho: \pi_1( X, x_0) \to K_{G/Z_G}$.
  The canonical homomorphism $\Aut( E) \to \Aut( E/Z_G)$ has kernel $Z_G$, so we get an embedding
  \begin{equation*}
    \Aut( E)/Z_G \hookrightarrow \Aut( E/Z_G) = C_{G/Z_G}( \rho) \subseteq G/Z_G
  \end{equation*}
  using also Lemma \ref{lemma:Aut}. Let $Z_{\varphi} \subseteq G/Z_G$ denote the nontrivial finite subgroup generated by
  the image of $\varphi \in \Aut( E)$ under this map.

  We choose a maximal torus $T \subseteq G$ such that $T/Z_G$ contains $Z_{\varphi}$, and decompose $\mathfrak g$
  into $Z_{\varphi}$-eigenspaces ${\mathfrak g}_{\chi}$ with respect to characters $\chi: Z_{\varphi} \to \Gm$.
  The fundamental group $\pi_1( X, x_0)$ also acts on $\mathfrak g$, via the composition
  \begin{equation*}
    \pi_1( X, x_0) \longto K_{G/Z_G} \subseteq G/Z_G \longto \Aut( \mathfrak g)
  \end{equation*}
  of $\rho$ with the adjoint action of $G/Z_G$ on $\mathfrak g$, and we have
  \begin{equation*}
    \ad( E) = \Xtilde \times^{\rho} {\mathfrak g}.
  \end{equation*}
  Because this action of $\pi_1( X, x_0)$ on $\mathfrak g$ commutes with the adjoint action of
  $Z_{\varphi} \subseteq G/Z_G$ on $\mathfrak g$, we get a vector bundle decomposition
  \begin{equation*}
    \ad( E) = \bigoplus_{\chi: Z_{\varphi} \to \Gm} \ad( E)_{\chi} \quad\text{with}\quad
      \ad( E)_{\chi} = \Xtilde \times^{\rho} {\mathfrak g}_{\chi}.
  \end{equation*}
  Since $Z_{\varphi}$ is nontrivial, there is a root $\alpha: T/Z_G \to \Gm$ of $G$ which is nontrivial on $Z_{\varphi}$.
  Then $Z_{\varphi}$ acts on the root spaces ${\mathfrak g}_{\alpha}$ and ${\mathfrak g}_{-\alpha}$ by nontrivial characters.
  Hence the eigenspace ${\mathfrak g}_{\chi_0}$ for the trivial character $\chi_0: Z_{\varphi} \to \Gm$
  has codimension $\geq 2$ in ${\mathfrak g}$. Thus the subbundle
  \begin{equation*}
    \ad( E)_{\chi_0} \subseteq \ad( E)
  \end{equation*}
  is a direct summand of corank $\geq 2$. Because $E$ is stable, $\ad( E)$ is semistable of degree $0$.
  It follows that $\ad( E)_{\chi_0}$ also has degree $0$. Thus Riemann-Roch implies that the linear subspace
  \begin{equation*}
    V^{\varphi} = H^1( X, \ad( E)_{\chi_0}) \subseteq V = H^1( X, \ad( E))
  \end{equation*}
  has codimension $\geq 2( g_X - 1) \geq 4$. So Lemma \ref{lemma:finitequotient} applies.  
\end{proof}

\section{The strictly semistable locus} \label{sec:str_semistab}
Let $G$ still be a reductive connected linear algebraic group over $\bbC$. Let $H$ be a Levi subgroup of a maximal parabolic subgroup $P \subsetneq G$.
Let $K_H \subseteq H$ be a maximal compact subgroup. We choose a Borel subgroup $B \subseteq P$ of $G$, and a maximal torus $T \subseteq B \cap H$ of $G$. Let
\begin{equation*}
  \alpha_1, \ldots, \alpha_r: T \longto \Gm \qquad\text{and}\qquad \alpha_1^{\vee}, \ldots, \alpha_r^{\vee}: \Gm \longto T
\end{equation*}
be the simple roots and coroots of $G$.
The quotient of $\Hom( \Gm, T)$ modulo its subgroup $\bbZ \cdot \alpha_1^{\vee} + \cdots + \bbZ \cdot \alpha_r^{\vee}$ is isomorphic to $\pi_1( G)$.

The Dynkin diagram of $H$ is given by removing one simple root $\alpha_i$ from the Dynkin diagram of $G$. The centers of $G$ and $H$ are
\begin{equation*}
  Z_G = \bigcap_j \ker( \alpha_j) \qquad\text{and}\qquad Z_H = \bigcap_{j \neq i} \ker( \alpha_j);
\end{equation*}
cf. \cite[(13.4)]{Hu}. It follows that $\alpha_i$ induces an isomorphism
\begin{equation*}
  Z_H/Z_G \longto[ \sim] \Gm.
\end{equation*}
The simple coroots of $H$ are the $\alpha_j^{\vee}$ with $j \neq i$, so the sequence
\begin{equation*}
  0 \longto \bbZ \longto[ \cdot \alpha_i^{\vee}] \pi_1( H) \longto \pi_1( G) \longto 0
\end{equation*}
is exact. Hence the induced map $\pi_1( [H, H]) \to \pi_1( [G, G])$ is injective.

Given a principal $H$--bundle $E$ over $X$, we denote by
\begin{equation*}
  E_G := E \times^H G = (E \times G)/H \longto X
\end{equation*}
the principal $G$--bundle obtained by extending the structure group of $E$ using the inclusion map $H \hookrightarrow G$.
If $E$ is semistable and its topological type $e \in \pi_1( H)$ is torsion, then $E_G$ is semistable \cite[Lemma 3.5.11]{Ra_Proc}.
Sending $[E]$ to $[E_G]$ thus defines a morphism of projective schemes
\begin{equation*}
  i_H: \calM_H' \longto \calM_G'.
\end{equation*}

The automorphism group of $H$ acts on $\calM_H$, by extension of the structure group.
An inner automorphism of $H$ acts trivially on $\calM_H$, so the group $\Out( H)$ of outer automorphisms also acts on $\calM_H$.

We denote the normalizer of $H$ in $G$ by $N_G( H)$. It acts on $H$ by conjugation in $G$. Let $\Gamma_H$ denote the image of $N_G( H)$ in $\Out( H)$.
The induced action of $\Gamma_H$ on $Z_H/Z_G \cong \Gm$ is effective,
since the centralizer of $Z_H/Z_G$ in $G/Z_G$ is known to be $H/Z_G$. This yields an embedding
\begin{equation*}
  \Gamma_H \hookrightarrow \Aut( \Gm) = \{ \pm 1\}.
\end{equation*}
The morphism $i_H$ is $\Gamma_H$-invariant, so it descends to a morphism
\begin{equation*}
  \ibar_H: \calM_H'/\Gamma_H \longto \calM_G'.
\end{equation*}
Note that $\Gamma_H$ acts trivially on $\pi_1( [H, H])$, because $N_G( H)$ acts trivially on $\pi_1( [G, G])$.
Hence $\Gamma_H$ acts on each component $\calM_H^e \subseteq \calM_H'$.
\begin{proposition} \label{normalization}
  Assume that $G$ is semisimple. Then the normalization of the strictly semistable locus
  \begin{equation*}
    \calM_G \setminus \calM_G^{\stab} \subset \calM_G
  \end{equation*}
  is given by the disjoint union
  \begin{equation*}
    \ibar := \amalg \ibar_H: \calM_H'/\Gamma_H \longto \calM_G' = \calM_G
  \end{equation*}
  over a set of representatives $H$ for conjugacy classes of Levi subgroups of maximal parabolic subgroups in $G$.
\end{proposition}
\begin{proof}
  As $\calM_H'$ is normal, its quotient $\calM_H'/\Gamma_H$ is also normal, due to \cite[p.~5]{Mu}.
  The image of $\ibar$ is contained in $\calM_G \setminus \calM_G^{\stab}$ by construction.
  Conversely, every closed point in $\calM_G \setminus \calM_G^{\stab}$ has the form
  \begin{equation*}
    [E_G] \in \calM_G \setminus \calM_G^{\stab}
  \end{equation*}
  for some semistable principal $H$--bundle $E$ over $X$ such that $E_G$ is also semistable. Let $e \in \pi_1( H)$ be the topological type of $E$. Let
  \begin{equation} \label{eq:chi_H}
    \chi_H: H \longto \Gm
  \end{equation}
  denote the generator of $\Hom( H, \Gm) \cong \bbZ$ which is dominant for $P$.
  If $\langle \chi_H, e \rangle > 0$, then the reduction $E_P := E \times^H P$ of $E_G$ to $P$ would violate semistability by \cite[Lemma 2.1]{Ra_Ann};
  if $\langle \chi_H, e \rangle < 0$, then the reduction $E_Q$ of $E_G$ to the opposite parabolic subgroup $H \subseteq Q \subseteq G$ would do so.
  These contradictions show $\langle \chi_H, e \rangle = 0$, which means $e \in \pi_1( [H, H])$.
  This proves that the image of $\ibar$ is equal to $\calM_G \setminus \calM_G^{\stab}$.

  The projective morphism $i_H$ is also affine due to Ramanathan's lemma \cite[Lemma 4.8.1, Remark 4.8.2 and Lemma 5.1]{Ra_Proc}.
  Hence the morphisms $i_H$, $\ibar_H$ and $\ibar$ are all finite.
  It remains to prove that $\ibar$ is generically injective. For that, we use Lemma \ref{lemma:dense}.

  Let $H'$ be a Levi subgroup of another maximal parabolic subgroup $P' \subseteq G$. Let $K_{H'} \subseteq H'$ again
  be a maximal compact subgroup. Suppose that two points
  \begin{equation*}
    [E_{\rho}] \in \calM_H' \qquad\text{and}\qquad [E_{\rho'}] \in \calM_{H'}'
  \end{equation*}
  have the same image in $\calM_G$ for some homomorphisms
  \begin{equation*}
    \rho: \pi_1( X, x_0) \longto K_H \qquad\text{and}\qquad \rho': \pi_1( X, x_0) \to K_{H'}.
  \end{equation*}
  Then there is an element $g \in G$ such that
  \begin{equation*}
    g \rho(\omega) g^{-1} = \rho'(\omega)
  \end{equation*}
  holds for all $\omega \in \pi_1( X, x_0)$. Now assume moreover that the image of $\rho$ is dense in $K_H$.
  Then we can conclude
  \begin{equation*}
    g K_H g^{-1} \subseteq K_{H'},
  \end{equation*}
  and hence $g H g^{-1} \subseteq H'$, so
  \begin{equation*}
    g H g^{-1} = H'
  \end{equation*}
  by maximality. Thus $H$ and $H'$ are conjugate in $G$, so we may assume $H = H'$, and then we get $g \in N_G( H)$.
  It follows that the two points $[E_{\rho}]$ and $[E_{\rho'}]$ are indeed in the same $\Gamma_H$-orbit.
\end{proof}

We return to the general case where $G$ is reductive. Let
\begin{equation*}
  \Sigma_G \subseteq \calM_G^{\sing}
\end{equation*}
denote the set of all singular closed points $[E] \in \calM_G$ such that every Euclidean neighborhood
$[E] \in U \subseteq \calM_G$ contains an open neighborhood $[E] \in U' \subseteq U$ for which $U' \setminus \calM_G^{\sing}$ is connected and simply connected.
Proposition \ref{prop:fin_quot} shows $\Sigma_G \subseteq \calM_G \setminus \calM_G^{\stab}$.

\begin{proposition} \label{prop:Gm_quot}
  Let $E$ be a principal $H$--bundle over $X$ such that $E_G$ is polystable with $\Aut( E_G) = Z_H$.
  Then the point $[E_G] \in \calM_G$ is in $\Sigma_G$.
\end{proposition}
\begin{proof}
  Luna's \'{e}tale slice theorem \cite{Lu} implies that $\calM_G$ is near the point $[E_G]$ analytically isomorphic to the GIT--quotient
  \begin{equation*}
    H^1( X, \ad( E_G))/\!/\Aut( E_G)
  \end{equation*}
  near $0$. Here $\Aut( E_G) = Z_H$ acts on the vector space
  \begin{equation*}
    V := H^1( X, \ad( E_G)) = H^1( X, E \times^H {\mathfrak g})
  \end{equation*}
  via the adjoint action of $Z_H$ on $\mathfrak g$. Since $Z_G \subseteq Z_H$ acts trivially,
  the isomorphism $\alpha_i: Z_H/Z_G \to \Gm$ yields an action of $\Gm$ on $\mathfrak g$.
  Let ${\mathfrak g}_m \subseteq {\mathfrak g}$ denote the $\Gm$-eigenspace of weight $m \in \bbZ$.
  Since the adjoint action of $H$ on $\mathfrak g$ commutes with $\Gm$, we get a vector bundle decomposition
  \begin{equation*}
    \ad( E_G) = \bigoplus_{m \in \bbZ} \ad( E_G)_m \quad\text{with}\quad \ad( E_G)_m := E \times^H {\mathfrak g}_m.
  \end{equation*}
  As $E_G$ is semistable, $\ad( E_G)$ is semistable of degree $0$. Consequently, $\ad( E_G)_m$ also has degree $0$. Due to Riemann-Roch, the $\Gm$-eigenspace
  \begin{equation*}
    V_m = H^1( X, \ad( E_G)_m) \subseteq V
  \end{equation*}
  has dimension $\geq (g_X - 1) \dim {\mathfrak g}_m$. Using ${\mathfrak g}_{\pm \alpha_i} \subseteq {\mathfrak g}_{\pm 1}$, we conclude
  \begin{equation*}
    \dim V_{\pm 1} \geq g_X - 1 \geq 2.
  \end{equation*}
  Then Lemma \ref{lemma:Gm} completes the proof.
\end{proof}
\begin{corollary} \label{cor:str_semistab}
  The strictly semistable locus $\calM_G \setminus \calM_G^{\stab}$ is
  the Zariski closure of the subset $\Sigma_G \subseteq \calM_G$.
\end{corollary}
\begin{proof}
  The image $K_{H/Z_G} \subseteq H/Z_G$ of $K_H \subseteq H$ is again a maximal compact subgroup.
  Let $E$ be a principal $H$--bundle over $X$ such that the principal $H/Z_G$--bundle $E/Z_G$ comes from a homomorphism
  \begin{equation*}
    \rho: \pi_1(X,x_0) \longto K_{H/Z_G}.
  \end{equation*}
  Then the principal $G/Z_G$--bundle $E_G/Z_G$ comes from the composition
  \begin{equation*}
    \pi_1( X, x_0) \longto[ \rho] K_{H/Z_G} \hookrightarrow K_{G/Z_G}
  \end{equation*}
  where $K_{G/Z_G} \subseteq G/Z_G$ is an arbitrary maximal compact subgroup containing $K_{H/Z_G}$.
  Hence $E_G/Z_G$ is polystable, so $E_G$ is also polystable.

  If moreover the image of $\rho$ is dense in $K_{H/Z_G}$, then Lemma \ref{lemma:Aut} yields
  \begin{equation*} 
    \Aut( E_G/Z_G) = C_{G/Z_G}( \rho) = Z_H/Z_G
  \end{equation*}
  since the centralizer of $H/Z_G$ in $G/Z_G$ is known to be $Z_H/Z_G$. This implies that the two canonical embeddings
  \begin{equation*}
    Z_H \hookrightarrow \Aut( E_G) \qquad\text{and}\qquad \Aut( E_G)/Z_G \hookrightarrow \Aut( E_G/Z_G)
  \end{equation*}
  are both isomorphisms. Using Proposition \ref{prop:Gm_quot}, it follows that the point $[E_G] \in \calM_G$ is in $\Sigma_G$.
  This shows that $\Sigma_G \subseteq \calM_G$ contains the inverse image of every point in $\calM_{G/Z_G}$ that comes from 
  a homomorphism $\rho: \pi_1(X,x_0) \to K_{H/Z_G}$ with dense image. Such points are dense in
  \begin{equation*}
    \calM_{G/Z_G} \setminus \calM_{G/Z_G}^{\stab}
  \end{equation*}
  by Lemma \ref{lemma:dense} and Proposition \ref{normalization}, so $\Sigma_G$ is dense in $\calM_G \setminus \calM_G^{\stab}$.
\end{proof}
\begin{corollary} \label{cor:regstab}
  The smooth locus of $\calM_G$ consists precisely of the moduli points $[E] \in \calM_G$
  of regularly stable principal $G$--bundles $E$ over $X$.
\end{corollary}
\begin{proof}
  The singular locus $\calM_G^{\sing} \subseteq \calM_G$ is closed and contains $\Sigma_G$, so it contains
  the strictly semistable locus $\calM_G \setminus \calM_G^{\stab}$ due to the previous Corollary \ref{cor:str_semistab}.
  The rest follows from Proposition \ref{prop:fin_quot}.i.
\end{proof}

\section{Reconstructing the Riemann surface} \label{sec:proof}
In this section, we prove the following main result of this paper.
\begin{theorem} \label{thm:main}
Let $X, X'$ be compact connected Riemann surfaces of genus $g_X, g_{X'} \geq 3$.
Let $G, G'$ be nonabelian connected reductive linear algebraic groups over $\bbC$. 
Let $d \in \pi_1( G)$ and $d' \in \pi_1( G')$ be given.

If the smooth locus of $\calM_G^d( X)$ is isomorphic to the smooth locus of $\calM_{G'}^{d'}( X')$,
then $X$ is isomorphic to $X'$.
\end{theorem}
This theorem is proved here by reconstructing $X$ from the smooth locus of the variety $\calM_G^d( X)$.
Starting from a moduli space $\calM_{G'}^{d'}( X')$ with isomorphic smooth locus,
this reconstruction will automatically yield an isomorphic Riemann surface, thereby proving $X' \cong X$.

Thus it suffices to consider only one Riemann surface $X$.
We reconstruct $X$ from the smooth locus of $\calM_G^d = \calM_G^d( X)$, in several steps.

\subsection*{Step 1.}
Let $Z_G^0 \subseteq Z_G$ be the identity component. We put $\Gbar := G/Z_G^0$. The projection $G \twoheadrightarrow \Gbar$ induces a morphism
\begin{equation} \label{eq:projection}
  \pr: \calM_G^d \longto \calM_{\Gbar}^{\dbar}
\end{equation}
where $\dbar \in \pi_1( \Gbar)$ denotes the image of $d \in \pi_1( G)$.
\begin{lemma} \label{lemma:pr_*}
  We have $\pr_*( \calO_{ \calM_G^d}) = \calO_{ \calM_{\Gbar}^{\dbar}}$.
\end{lemma}
\begin{proof}
  The canonical embedding
  \begin{equation*}
    \calO_{ \calM_{\Gbar}^{\dbar}} \hookrightarrow \pr_*( \calO_{ \calM_G^d})
  \end{equation*}
  turns the latter into a coherent sheaf of algebras over the former. Let
  \begin{equation*}
    U \subseteq \calM_{\Gbar}^{\dbar}
  \end{equation*}
  be the regularly stable locus. This locus $U$ is known to be open and non-empty;
  one way to see this is to use Corollary \ref{cor:regstab}.
  Because $\calM_{\Gbar}^{\dbar}$ is normal, it suffices to prove the claim over $U$.

  Since $Z_G^0$ is central in $G$, the multiplication map
  \begin{equation*}
    Z_G^0 \times G \longto G
  \end{equation*}
  is a group homomorphism. It induces a morphism of projective varieties
  \begin{equation*}
    \calM^0_{Z_G^0} \times \calM_G^d \longto \calM_G^d
  \end{equation*}
  by extension of the structure group. Thus the abelian variety
  \begin{equation*}
    A := \calM^0_{Z_G^0}
  \end{equation*}
  acts on $\calM_G^d$. The map $\pr$ in \eqref{eq:projection} is $A$--invariant, and its restriction
  \begin{equation*}
    \pr: \pr^{-1}( U) \longto U
  \end{equation*}
  is a principal $A$--bundle. As $A$ is integral, the claim follows over $U$.
\end{proof}
Let $\omega_G$ denote the dualizing sheaf of $\calM_G^d$. By \cite[Proposition 2.2]{Pa} or \cite[Theorem 2.8]{KN}, $\omega_G$ is a line bundle on $\calM_G^d$.
Its restriction to the smooth locus $\calM_G^d \setminus \calM_G^{\sing}$ is the canonical line bundle $\det( \Omega^1)$.
\begin{corollary} \label{cor:glob_gen}
  Choose $n \gg 0$. Then the line bundle $\omega_G^{\otimes (-n)}$ on $\calM_G^d$ is globally generated,
  and the corresponding morphism
  \begin{equation*}
    \varphi_{\omega_G^{\otimes (-n)}}: \calM_G^d \longto \bbP^N
  \end{equation*}
  factors into the morphism $\pr$ in \eqref{eq:projection}, followed by a closed immersion.
\end{corollary}
\begin{proof}
 We have $\omega_G \cong \pr^*( \omega_{\Gbar})$ according to \cite[p. 525--527]{Pa}. Hence the projection formula gives an isomorphism
  \begin{equation*}
    \pr_* \big( \omega_G^{\otimes (-n)} \big) \cong \omega_{\Gbar}^{\otimes (-n)} \otimes \pr_*( \calO_{ \calM_G^d}) = \omega_{\Gbar}^{\otimes (-n)}
  \end{equation*}
  for any $n$. Taking global sections on both sides, we get
  \begin{equation*}
    H^0 \big( \calM_G^d, \omega_G^{\otimes (-n)} \big) = H^0 \big( \calM_{\Gbar}^{\dbar}, \omega_{\Gbar}^{\otimes (-n)} \big).
  \end{equation*}
  Now use that $\omega_{\Gbar}^{\otimes -1}$ is ample, according to \cite[Corollary 2.1]{Pa}.
\end{proof}
This shows how to reconstruct the morphism $\pr$ in \eqref{eq:projection}, and in particular its target
$\calM_{\Gbar}^{\dbar}$, from the variety $\calM_G^d$.

One can in fact reconstruct $\calM_{\Gbar}^{\dbar}$ from just the smooth locus of $\calM_G^d$, by taking the closure of its image under the map
\begin{equation*}
  \varphi_{\omega_G^{\otimes (-n)}}: \calM_G^d \setminus \calM_G^{\sing} \longto \bbP^N
\end{equation*}
given by the canonical line bundle $\omega_G$ on $\calM_G^d \setminus \calM_G^{\sing}$.
\begin{remark}
  Every holomorphic section of $\omega_G$ over $\calM_G^d \setminus \calM_G^{\sing}$ extends to $\calM_G^d$ by normality, and hence is algebraic.
  Thus it suffices to assume in Theorem \ref{thm:main} that the two smooth loci are biholomorphic.
\end{remark}
For the remaining steps, we can thus assume that $G$ is semisimple, $G \neq \{1\}$, and that we are given the projective variety $\calM_G^d$.

\subsection*{Step 2.} Corollary \ref{cor:str_semistab} characterizes the locus $\calM_G^d \setminus \calM_G^{\stab}$ in $\calM_G^d$.
We assume that it is non-empty, and consider an irreducible component of it. Its normalization is, due to Proposition \ref{normalization}, isomorphic to
\begin{equation*}
  \calM_H^e/\Gamma_H
\end{equation*}
for some Levi subgroup $H$ of a maximal parabolic subgroup $P \subsetneq G$, and some inverse image $e \in \pi_1([H, H])$ of $d$.

The group $\Gamma_H$ acts on the line bundle $\omega_H$ over $\calM_H^e$. Its square $\omega_H^{\otimes 2}$ descends to a line bundle $\omega_H^{\otimes 2}/\Gamma_H$
over $\calM_H^e/\Gamma_H$ according to Kempf's lemma \cite[Th\'{e}or\`{e}me 2.3]{DN}.
This line bundle is determined by the normal variety $\calM_H^e/\Gamma_H$ alone, due to the following lemma.
\begin{lemma}
  The line bundle $\omega_H^{\otimes 2}/\Gamma_H$ is over the smooth locus of $\calM_H^e/\Gamma_H$
  isomorphic to the square of the canonical line bundle $\det( \Omega^1)$.
\end{lemma}
\begin{proof}
  This is clear if $\Gamma_H$ is trivial. So we may assume $\Gamma_H \cong \{\pm 1\}$.
  The character $\chi_H: H \to \Gm$ in \eqref{eq:chi_H} induces a $\Gamma_H$--equivariant morphism
  \begin{equation*}
    (\chi_H)_*: \calM_H^e \longto J = J( X) := \Pic^0( X).
  \end{equation*}
  Since $\Gamma_H$ acts effectively on $Z_H/Z_G$, it also acts effectively on the isogenous torus $H/[H, H]$.
  Hence $-1 \in \Gamma_H$ acts on $J$ by the inverse map from the group structure on $J$. In particular, the fixed locus
  \begin{equation*}
    (\calM_H^e)^{\Gamma_H} \subseteq \calM_H^e
  \end{equation*}
  is contained in the inverse image of the $2$--division points in $J$, so its codimension is at least $g_X \geq 3$.
  Thus it suffices to check the lemma over $U/\Gamma_H$ for the smooth open subscheme
  \begin{equation*}
    U := \calM_H^e \setminus \big( (\calM_H^e)^{\sing} \cup (\calM_H^e)^{\Gamma_H} \big)
  \end{equation*}
  where $\Gamma_H$ acts freely. But here the claim follows simply from the fact that the canonical projection $U \twoheadrightarrow U/\Gamma_H$ is \'{e}tale.
\end{proof}

\subsection*{Step 3.} Replacing $G$ by $H$ in \eqref{eq:projection},
we get a canonical $\Gamma_H$-equivariant morphism $\pr: \calM_H^e \to \calM_{\Hbar}^{\ebar}$ with $\Hbar := H/Z_H^0$. It induces a morphism
\begin{equation} \label{eq:projection_mod_Gamma}
  \prbar: \calM_H^e/\Gamma_H \longto \calM_{\Hbar}^{\ebar}/\Gamma_H.
\end{equation}
\begin{lemma} \label{lemma:fibration} 
  There is a dense open subscheme $\Ubar \subseteq \calM_{\Hbar}^{\ebar}/\Gamma_H$ such that
  \begin{equation*}
    \prbar: (\prbar)^{-1}( \Ubar) \longto \Ubar
  \end{equation*}
  is an \'{e}tale-locally trivial fibration, with fiber either the Jacobian $J$
  or its quotient $J/\{\pm 1\}$ modulo the inverse map from its group structure.
\end{lemma}
\begin{proof}
  Let $U \subseteq \calM_{\Hbar}^{\ebar}$ be the regularly stable locus. The restriction
  \begin{equation} \label{eq:J_bundle}
    \pr: \pr^{-1}( U) \longto U
  \end{equation}
  is a principal $J$--bundle, as in the proof of Lemma \ref{lemma:pr_*}.

  Suppose for the moment that the action of $\Gamma_H$ on $\calM_{\Hbar}^{\ebar}$ is effective.
  Then $\Gamma_H$ acts freely on some dense open subscheme $U' \subseteq U$.
  The principal $J$--bundle $\pr$ in \eqref{eq:J_bundle} descends to a principal $J$--bundle
  \begin{equation*}
    \pr^{-1}( U')/\Gamma_H \longto U'/\Gamma_H.
  \end{equation*}
  Hence the lemma holds for $\Ubar := U'/\Gamma_H \subseteq \calM_{\Hbar}^{\ebar}/\Gamma_H$ in this case.

  It remains to treat the case where the action of $\Gamma_H$ on $\calM_{\Hbar}^{\ebar}$ is not effective.
  This means $\Gamma_H \cong \{\pm 1\}$, and that $\Gamma_H$ acts trivially on $\calM_{\Hbar}^{\ebar}$.

  We claim that in this case $\Ubar := U \subseteq \calM_{\Hbar}^{\ebar} = \calM_{\Hbar}^{\ebar}/\Gamma_H$ satisfies the conditions in the lemma.
  Replacing $U$ by an \'{e}tale covering and pulling back, we can choose a section $\sigma$ of the principal $J$--bundle
  $\pr$ in \eqref{eq:J_bundle}. Let $\gamma \in \Gamma_H$ be the nontrivial element. Then
  \begin{equation*}
    \gamma \cdot \sigma = \xi + \sigma
  \end{equation*}
  for some morphism $\xi: U \to J$. Refining the covering of $U$ if necessary, and using the divisibility of $J$,
  we may assume $\xi = 2 \xi'$. Then
  \begin{equation*}
    \gamma \cdot (\xi' + \sigma) = \gamma \cdot \xi' + \gamma \cdot \sigma = - \xi' + (\xi + \sigma) = \xi' + \sigma,
  \end{equation*}
  since $\Gamma_H$ acts by the inverse map on $J$. So $\xi' + \sigma$ is $\Gamma_H$-invariant.

  This shows that the principal $J$--bundle $\pr$ in \eqref{eq:J_bundle} admits
  \'{e}tale-locally $\Gamma_H$-equivariant trivializations. Hence it descends to
  an \'{e}tale-locally trivial fibration with fiber $J/\Gamma_H = J/\{\pm 1\}$ over $U$.
\end{proof}
\begin{remark}
  The second case, where the general fiber of $\prbar$ is $J/\{\pm 1\}$, actually occurs.
  An example is $G = \Sp_{2n + 2}$ and $H \cong \Gm \times \Sp_{2n}$.
  Here $\Gamma_H \cong \{\pm 1\}$ is nontrivial, but its action on $\Hbar = \Sp_{2n}$ is trivial.
\end{remark}
\begin{corollary}
  Choose $n \gg 0$. Then the line bundle $(\omega_H^{\otimes 2}/\Gamma_H)^{\otimes( -n)}$ on $\calM_H^e/\Gamma_H$ is globally generated,
  and the corresponding morphism
  \begin{equation*}
    \varphi_{(\omega_H^{\otimes 2}/\Gamma_H)^{\otimes( -n)}}: \calM_H^e/\Gamma_H \longto \bbP^N
  \end{equation*}
  factors into the morphism $\prbar$ in \eqref{eq:projection_mod_Gamma}, followed by a closed immersion.
\end{corollary}
\begin{proof}
  Using Kempf's lemma, the anti-ample line bundle $\omega_{\Hbar}^{\otimes 2}$ on $\calM_{\Hbar}^{\ebar}$ descends
  to an anti-ample line bundle $\omega_{\Hbar}^{\otimes 2}/\Gamma_H$ on $\calM_{\Hbar}^{\ebar}/\Gamma_H$ with
  \begin{equation*}
    \prbar^*( \omega_{\Hbar}^{\otimes 2}/\Gamma_H) \cong \omega_H^{\otimes 2}/\Gamma_H.
  \end{equation*}
  The corollary can now be proved exactly as Corollary \ref{cor:glob_gen},
  since Lemma \ref{lemma:fibration} and the normality of $\calM_{\Hbar}^{\ebar}/\Gamma_H$ together imply
  \begin{equation*}
    \prbar_*( \calO_{\calM_H^e/\Gamma_H}) = \calO_{\calM_{\Hbar}^{\ebar}/\Gamma_H}. \qedhere
  \end{equation*}
\end{proof}
This shows how to reconstruct the morphism $\prbar$ in \eqref{eq:projection_mod_Gamma}, and in particular its 
general fiber $J( X)$ or $J( X)/\{ \pm 1\}$, from the variety $\calM_H^e/\Gamma_H$.

It is known that $J( X)$ can be reconstructed from the Kummer variety $J( X)/\{ \pm 1\}$,
for example as its integral closure in the field of meromorphic functions on the two-sheeted cover
of its smooth locus given by the unique maximal torsionfree subgroup in 
its fundamental group.

\subsection*{Step 4.} Starting from the variety $\calM_G^d$, we have now reconstructed the Jacobian $J = J( X)$, together with a morphism $j_E: J \to \calM_G^d$
given by a general point $[E] \in \calM_H^e$. This morphism is the composition
\begin{equation*}
  j_E: J \cong \calM_{Z^0_H}^0 \longto[ \cdot E] \calM_H^e \longto[i_H] \calM_G^d.
\end{equation*}
\begin{lemma}
  The first Chern class $c_1( j_E^* \, \omega_G) \in H^2( J, \bbZ)$ is a negative multiple of the  canonical principal polarization $\Theta$ on $J$.
\end{lemma}
\begin{proof}
  According to \cite[Proposition 2.2]{Pa}, $\omega_G$ is the pullback along
  \begin{equation*}
    \ad: \calM_G^d \longto \calM_{\SL( {\mathfrak g})}
  \end{equation*}
  of the determinant of cohomology line bundle $\Ldet$ on $\calM_{\SL( {\mathfrak g})}$.
  Points in $\calM_{\SL( {\mathfrak g})}$ correspond to polystable vector bundles $V$ over $X$ with trivial determinant;
  we recall that the fiber of $\Ldet$ over such a point is the determinant of the cohomology of $V$, i.\,e.\ the vector space
  \begin{equation*}
    \Lambda^{\rm top} H^0(X, V) \otimes \Lambda^{\rm top} H^1(X, V)^*
  \end{equation*}
  of dimension one. Let ${\mathfrak g}_m \subseteq {\mathfrak g}$ denote the eigenspace where $\Gm \cong Z_H^0$ acts with weight $m \in \bbZ$.
  The line bundle $j_E^* \, \omega_G$ is the pullback of $\Ldet$ along the morphism
  \begin{equation*}
    J \longto \calM_{\SL( {\mathfrak g})}
  \end{equation*}
  that sends a line bundle $L$ of degree $0$ over $X$ to the vector bundle
  \begin{equation*}
    \bigoplus_{m \in \bbZ} L^{\otimes m} \otimes ( E \times^H {\mathfrak g}_m)
  \end{equation*}
  with trivial determinant over $X$. Computing the determinant of its cohomology, we conclude that
  \begin{equation*}
    c_1( \calL) = \sum_{m \in \bbZ} - m^2 \dim {\mathfrak g}_m \cdot \Theta \in H^2( J, \bbZ).
  \end{equation*}
  This is a strictly negative multiple of $\Theta$, since ${\mathfrak g}_m \neq 0$ for some $m \neq 0$,
  because the action of $\Gm \cong Z_H^0$ on $\mathfrak g$ is not trivial.
\end{proof}
Hence we can also reconstruct the principal polarization $\Theta$ on $J( X)$; then the usual Torelli theorem gives back the Riemann surface $X$.

\subsection*{Step 5.} It remains to treat the case $\calM_G^d \setminus \calM_G^{\stab} = \emptyset$.
According to \cite[Proposition 7.8]{Ra_Ann}, this can only happen if
\begin{equation*}
  G/Z_G \cong \PGL_{n_1} \times \cdots \times \PGL_{n_r}.
\end{equation*}
In particular, the universal covering of $G$ is a homomorphism
\begin{equation} \label{eq:univ_cover}
  \SL_{n_1} \times \cdots \times \SL_{n_r} \longto G.
\end{equation}
Its kernel is a finite subgroup $\mu$ in the product of the centers $\mu_{n_i} \subseteq \SL_{n_i}$.

For each $i$, we lift the image of $d \in \pi_1( G)$ in $\pi_1( \PGL_{n_i}) \cong \bbZ/n_i \bbZ$ to an integer $d_i \in \bbZ$.
Then the emptiness of the strictly semistable locus also implies that $d_i$ is coprime to $n_i$ for all $i$.
Let $L_i$ be a line bundle of degree $d_i$ on $X$, and let $\calM_{n_i, L_i}$ denote the coarse moduli space of semistable
vector bundles $E$ over $X$ with rank $n_i$ and determinant $L_i$. The homomorphism \eqref{eq:univ_cover} induces a finite
surjective morphism
\begin{equation} \label{eq:projection_coprime}
  \pr: \calM_{n_1, L_1} \times \cdots \times \calM_{n_r, L_r} \longto \calM_G^d.
\end{equation}
The smooth locus in $\calM_G^d$ coincides with the regularly stable locus by Corollary \ref{cor:str_semistab};
we denote it by $U \subseteq \calM_G^d$.
The abelian group of principal $\mu$--bundles over $X$ is isomorphic to $\mu^{2 g_X}$.
This group acts on the morphism \eqref{eq:projection_coprime}. This turns the restriction
\begin{equation*}
  \pr: \pr^{-1}( U) \longto U
\end{equation*}
into an unramified Galois covering with group $\mu^{2 g_X}$.
Since $\calM_{n_i, L_i}$ is a unirational smooth projective variety, it is in particular simply connected \cite{Se}.
Thus $\pr^{-1}( U)$ is also simply connected, as we are removing a closed subset of complex codimension $\geq 2$.
This shows that $\pr^{-1}( U)$ is the universal covering of $U$. Taking the integral closure of $\calM_G^d$ in the
field of meromorphic functions on $\pr^{-1}( U)$, we can thus recover the morphism \eqref{eq:projection_coprime}
from the variety $\calM_G^d$ in this case.

Because the Picard group of $\calM_{n_i, L_i}$ is isomorphic to $\bbZ$, the extremal rays in the nef cone of the variety
\begin{equation*}
  \calM_{n_1, L_1} \times \cdots \times \calM_{n_r, L_r}
\end{equation*}
correspond to the factors $\calM_{n_i, L_i}$.
Thus we can even reconstruct these factors $\calM_{n_i, L_i}$ from the variety $\calM_G^d$.

The projective variety $\calM_{n_i, L_i}$ is known to determine the Riemann surface $X$ up to isomorphy; 
cf. \cite[Theorem 1]{T2} or \cite[Theorem 3]{NR_Def}.
Thus we can reconstruct $X$ from $\calM_G^d( X)$ in this case as well.

This completes the proof of Theorem \ref{thm:main}. \qed

\end{document}